\documentclass{amsart}

\usepackage{amsfonts,amssymb,verbatim,amsmath,amsthm,latexsym,textcomp,amscd}
\usepackage{latexsym,amsfonts,amssymb,epsfig,verbatim}
\usepackage{amsmath,amsthm,amssymb,latexsym,graphics,textcomp}
\usepackage{graphicx}
\usepackage{color}
\usepackage{url}

\begin{document}

\newtheorem{theorem}{Theorem}[section]
\newtheorem{prop}[theorem]{Proposition}
\newtheorem{lemma}[theorem]{Lemma}
\newtheorem{cor}[theorem]{Corollary}
\newtheorem{defn}[theorem]{Definition}
\newtheorem{conj}[theorem]{Conjecture}
\newtheorem{claim}[theorem]{Claim}
\newtheorem{example}[theorem]{Example}
\newtheorem{rem}[theorem]{Remark}
\newtheorem{rmk}[theorem]{Remark}
\newtheorem{obs}[theorem]{Observation}

\newcommand{\A}{\rightarrow}
\newcommand{\C}{\mathcal C}
\newcommand\AAA{{\mathcal A}}
\newcommand\BB{{\mathcal B}}
\newcommand\DD{{\mathcal D}}
\newcommand\EE{{\mathcal E}}
\newcommand\FF{{\mathcal F}}
\newcommand\GG{{\mathcal G}}
\newcommand\HH{{\mathcal H}}
\newcommand\I{{\stackrel{\rightarrow}{i}}}
\newcommand\J{{\stackrel{\rightarrow}{j}}}
\newcommand\K{{\stackrel{\rightarrow}{k}}}
\newcommand\LL{{\mathcal L}}
\newcommand\MM{{\mathcal M}}
\newcommand\NN{{\mathbb N}}
\newcommand\OO{{\mathcal O}}
\newcommand\PP{{\mathcal P}}
\newcommand\QQ{{\mathcal Q}}
\newcommand\RR{{\mathcal R}}
\newcommand\SSS{{\mathcal S}}
\newcommand\TT{{\mathcal T}}
\newcommand\UU{{\mathcal U}}
\newcommand\VV{{\mathcal V}}
\newcommand\WW{{\mathcal W}}
\newcommand\XX{{\mathcal X}}
\newcommand\YY{{\mathcal Y}}
\newcommand\ZZ{{\mathbb Z}}
\newcommand\hhat{\widehat}
\newcommand\vfn{\stackrel{\A}{r}(t)}
\newcommand\dervf{\frac{d\stackrel{\A}{r}}{dt}}
\newcommand\der{\frac{d}{dt}}
\newcommand\vfncomp{f(t)\I+g(t)\J+h(t)\K}
\newcommand\ds{\sqrt{f^{'}(t)^2+g^{'}(t)^2+h^{'}(t)^2}dt}
\newcommand\rvec{\stackrel{\A}{r}}
\newcommand\velo{\frac{d\stackrel{\A}{r}}{dt}}
\newcommand\speed{|\velo|}
\newcommand\velpri{\rvec \,^{'}}

\title[Bounded packing in CAT(0) spaces and Gromov hyperbolic spaces]{ Existence and non-existence of bounded packing in CAT(0) spaces and Gromov hyperbolic spaces}

\author{Pranab Sardar}
\address{Chennai Mathematical Institute, India
\newline
Email: pranab@cmi.ac.in}

\date{\today}

\begin{abstract}
The main result of this paper is that given a group $G$ acting geometrically by isometries on a CAT(0)
space $X$
and a cyclic subgroup $H$ of $G$ generated by a rank-1 isometry of $X$, $H$ has bounded packing in $G$.
We give two proofs of this result. The first one is by a clever argument of Mj(Lemma 3.3 of \cite{mahan-CC})
and the characterization of rank-$1$ isometries by Hamenstadt(\cite{hamen-rank1}). The second proof
follows directly from some results of Dahmani, Guirardel and Osin(\cite{dah-gui-osin})
and Sisto(\cite{sisto}). Then using Mihailova's construction, we
show the existence of a finitely generated subgroup of the direct product of two free groups
$\mathbb F_2\times \mathbb F_2$
without the bounded packing property answering a question of Hruska-Wise(\cite{hruska-wise}). We also
prove the existence of finitely presented subgroups of CAT(0) groups without bounded packing using
Wise's {\em modified Rip's construction}(\cite{wise-rips1}) and the {\bf 1-2-3} theorem
of Baumslag, Bridson, Miller and Short(\cite{bbms1}).
\end{abstract}

\maketitle

\section{Introduction}

Bounded packing was defined by Hruska and Wise (\cite{hruska-wise}) motivated by the concept of width
of subgroups due to \cite{GMRS}. Hruska and Wise proved bounded packing of quasi-convex subgroups of 
Gromov hyperbolic groups and relatively quasi-convex subgroups of relatively hyperbolic groups under
mild and natural restrictions. Bounded packing proves to be useful for two reasons. One is that
failure of bounded packing of $H$ in $G$- which is a geometric condition, implies that $H$ is not
separable in $G$- which is an algebraic property. This result is due to Yang(\cite{yang}). On the other
hand, as commented by Hruska-Wise  in \cite{hruska-wise}(Corollary $3.1$) the following result
follows from the work of Sageev(\cite{sageev-codim1}):

\smallskip
{\bf Theorem:} {\em Suppose $H$ is a finitely generated codimension–$1$ subgroup
of a finitely generated group $G$. If $H$ has bounded packing in $G$, then the corresponding
CAT(0) cube complex $C$ is finite dimensional.}

\smallskip
However, constructing new examples of groups with subgroups with (or without) bounded packing seems
difficult although it is known to hold for quasi-convex subgroups of hyperbolic groups, all subgroups of
polycyclic groups and nilpotent groups and so on. (See Example 2.19, Example 2.22 of \cite{hruska-wise},
and the main results of \cite{yang} and \cite{bddpacking} for more examples.)
Moreover, many natural questions about bounded packing remain unanswered.For example in 
\cite{hruska-wise} we have:

\smallskip
{\bf Problem}(Problem 2.24 of \cite{hruska-wise}): {\em Give an example of a cyclic subgroup $Z$ of a
finitely generated group $G$ such that $Z$ does not have bounded packing in $G$. }

\smallskip
Recently Wise and Woodhouse (\cite{wise-wood}) have shown that
abelian subgroups of groups acting geometrically on CAT(0) cube complexes have bounded packing.
However, the following special case of the above problem is still not answered.

\smallskip
{\bf Problem} (Problem 1.1(1) of \cite{wise-wood}): {\em
Let $G$ act properly and cocompactly on a CAT(0) space. Does each [cyclic]
abelian subgroup $A$ of  $G$ have bounded packing?}

\smallskip
Towards this direction we show that if $G$ is a group acting geometrically a CAT(0) space $X$ and
$H$ is a cyclic subgroup generated by a rank-$1$ isometry then $H$ has bounded packing in $G$.
However, we prove that there are finitely
generated subgroups of $\mathbb F_2\times \mathbb F_2$ without bounded packing answering Question $2.25$ of 
\cite{hruska-wise}. Also we show that there are finitely presented subgroups of CAT(0) groups without
bounded packing using Wise's {\em modified Rip's construction}(\cite{wise-rips1}) and the 
{\bf 1-2-3} theorem of Baumslag, Bridson, Miller and Short(\cite{bbms1}).

{\bf Acknowledgements:} I would like to thank Misha Kapovich for many many useful discussions.
In particular,
he told me to look at Mihailova's constructions to find examples of subgroups of $F_2\times F_2$
without bounded packing and to look at Hamnestadt's paper \cite{hamen-rank1}.
I also would like to thank Mahan Mj for inspiration and helpful email exchanges and for pointing me
to the Lemma $3.3$ of \cite{mahan-CC}. Most part of this work was done during my term in the
University of California, Davis. I would like to thank University of California, Davis and
Chennai Mathematical Institute, India for hospitality and financial support.

\section{Bounded Packing}
\begin{defn}(\cite{hruska-wise})
Suppose $G$ is a countable group with a proper, left invariant metric $d$. Let $H\leq  G$ be a subgroup.
We say that $H$ has bounded packing in $G$ (with respect to $d$) if for any $D>0$ there is a number
$n=n(G,H,D)$ such that given any collection of left cosets $I\subset G/H$ of $H$ in $G$ such
that $d(g_1H, g_2H)\leq D$ for all $g_1H,g_2H \in I$, we have $|I|\leq n$.
\end{defn}

By Lemma 2.2 of \cite{hruska-wise} we know that bounded packing for a subgroup $H$ of a countable group $G$ 
is independent of the choice of the particular left invariant proper metric on $G$. On the other hand, by the work of Higman,
Neumann and Neumann(\cite{hnn}) we know that any countable group can be embedded in a $2$-generated group. 
Since every finitely generated group admits a proper, invariant metric- namely a word metric, it follows that any countable 
group admits a left invariant proper metric. Therefore, bounded packing of subgroups makes sense for subgroups of all countable groups. 

The following lemma is a summary of some basic results about bounded packing which are proved
by Hruska-Wise in the second section of \cite{hruska-wise}. We will make repeated use of it later.

\begin{lemma}\label{lemma1} (\cite{hruska-wise}) Let $G$ be a countable group with a proper left invariant metric.
Then the following are true.
\begin{enumerate}
\item Every finite subgroup of $G$ has bounded packing in $G$. 
\item Every finite index subgroup of $G$ has bounded packing in $G$.
\item For any sequence of subgroups $K\subset H\subset G$ if $K$ has bounded packing in $H$ and 
$H$ has bounded packing in $G$, then $K$ has bounded packing in $G$.
\item For any sequence of subgroups $K\subset H\subset G$ if $[H:K]<\infty$
then $K$ has bounded packing in $G$ if and only if $H$ has bounded packing in $G$.
\item If $H,K$ are two subgroup of $G$ both of which have bounded packing in $G$ then so does their intersection
$H\cap K$.
\item Suppose $H, \,K \leq G$ and $[G : K] < \infty$ . Then $H\cap K$ has bounded packing in $K$
if and only if $H$ has bounded packing in $G$.
\item Suppose we have a short exact sequence of groups $1\rightarrow N\rightarrow G\stackrel{\pi}{\rightarrow} Q\rightarrow 1$.
Let $Q_1$ be a subgroup of $Q$. Then $Q_1$ has bounded packing in $Q$ if and only if $\pi^{-1}(Q_1)$ has the bounded packing in $G$.

In particular, normal subgroups of {\em any} countable group $G$ have bounded packing in $G$.
\end{enumerate}
\end{lemma} 

\section{The main result}
The goal of this section is to prove the following theorem:

\begin{theorem}\label{th1}
Suppose $G$ is a finitely generated group acting by homeomorphism on a compact space $X$.
Suppose $H\leq G$ and $A\subset X$ invariant under $H$ such that the following conditions hold:
\begin{enumerate}
\item {\bf Limit set property:} For any $x\in X$ and any infinite sequence of distinct elements $\{h_n\}$ in $H$ there is
a subsequence $\{h_{n_k}\}$ such that $\lim_{k\rightarrow \infty} h_{n_k}.x $ exists and it is in $A$.
\item {\bf Dynamical quasi-convexity:} For any sequence of distinct cosets $\{g_iH\}$ there is a subsequence $\{g_{i_n}\}$ such
that $\lim_{n\rightarrow \infty} g_{i_n}.A$ is a single point in $X$.
\item For any pair of distinct cosets $g_1H, g_2H$ we have
$g_1A\cap g_2A=\emptyset$.
\end{enumerate}

Then $H$ has bounded packing in $G$.
\end{theorem}

The proof is an adaptation of the arguments in the proof of Lemma 3.3 of \cite{mahan-CC}.
For convenience it is broken into the following two lemmas. We assume that the first and the
second  conditions of the theorem hold.

\begin{lemma}
Suppose $C_i$ is an infinite sequence of sets of cosets of $H$ in $G$ such that for all $i$ and 
$xH,yH\in C_i$, $d(xH,yH)\leq D$ and $|C_i|\rightarrow \infty$. Then there is an infinite sequence
of distinct cosets $\{g_nH\}$ such that $d(1, g_1H)\leq D$, $d(1, g_2H)\leq D$, 
$d(1,g_nH) \rightarrow \infty$ and $d(g_iH,g_jH)\leq D$ for  $i=1,2$ and $j\geq 3$.
\end{lemma}

$Proof:$ For all $i$ we can find an element $x_i\in G$ such that under left multiplication by $x_i$
two cosets of $C_i$ intersect $B(1;D)$. Thus we get a sequence of sets of cosets of $H$ such that
each one of them has two members close to identity. Now, since the metric on $G$ is proper only
finitely many cosets of $H$ can intersect $B(1;D)$. Therefore, after passing to a subsequence we
may assume that each of $C^{'}_i:= x_iC_i$ contains two fixed cosets intersecting $B(1;D)$. Call
these cosets $g_1H,g_2H$. Now, since $|C_i|\rightarrow \infty$ we can pick one coset $g_nH$ from
$C^{'}_n$ such that $d(1,g_nH)\rightarrow \infty$ because of the properness of the metric on $G$.
Also, by construction $d(g_iH,g_jH)\leq D$ for $i=1,2$ and $j\geq 3$. $\Box$

\begin{lemma}
Suppose we have an an infinite sequence of distinct cosets $\{g_nH\}$ such that
$d(1, g_1H)\leq D$, $d(1,g_2H)\leq D$ and $d(g_iH,g_jH)\leq D$ for
$i=1,2$ and $j\geq 3$. Then $g_1A\cap g_2A\neq \emptyset$.
\end{lemma}

$Proof:$
We have a sequence of distinct cosets $\{g_nH\}$ such that there are points 
$v_{ij}\in g_iH, w_{ij} \in g_jH$ for $i= 1,2$, and $j\geq 3$ such that 
$d(v_{ij}, w_{ij})\leq D$. Therefore, we can write $w_{ij}=v_{ij}.z_{ij}$
where $d(1,z_{ij})\leq D$. If necessary, by passing to a subsequence, we may assume that 
$g_n.A$ converges to $p\in X$. Now fix a point $q\in A$ and $i$, $1\leq i \leq 2$. Without loss of generality we can assume that $v_{ij}.q$ is a convergent sequence.
Passing to a subsequence we can assume that $z_{ij}=z$. Write $v_{ij}=g_i.h_{ij}$ where $h_{ij}\in H$.
Now consider the sequence $w_{ij}.q=(g_i h_{ij}z).q=g_i(h_{ij}.(z.q))$. By the limit set property of $H$
we know that this sequence has a limit in $g_iA$. Thus $p\in g_iA$ for $i=1,2$. $\Box$

{\em Proof of the theorem:} Suppose $H$ does not have bounded packing in $G$. Then
there is a constant $D>0$ and an infinite sequence of sets $C_i$ of cosets of $H$ in
$G$ such that $|C_i|\rightarrow \infty$ and for all $i$ and $xH,yH\in C_i$ we have
$d(xH,yH)\leq D$.

Now, by the above lemmas then we can find two distinct cosets $g_1H, g_2H$ such that
$g_1A\cap g_2A\neq \emptyset$. This contradicts the third hypothesis of the theorem and
we are done. $\Box$

In \cite{hruska-wise} the authors raises the question
if relatively quasi-convex subgroups of relatively hyperbolic groups have finite width. Since
finite width is a consequence of bounded packing (see the proof of Corollary 8.11 in \cite{hruska-wise}),
we obtain a partial answer to this question in the form of the following corollary.

\begin{defn}
Suppose $G$ is a relatively hyperbolic group and $H$ is a subgroup. Let $\partial G$ be the
Bowditch boundary (see \cite{bowditch-relhyp} for detail) of $G$ and let $A$ be the limit set of $H$
in $\partial G$. We will call $H$ a {\em dynamically malnormal} subgroup of $G$ if for all
$g\in G\setminus H$ we have $gA\cap A=\emptyset$.
\end{defn}

Note that the same notion can be defined for a hyperbolic group $G$ taking the action of it
on the Gromov boundary. However, for a quasi convex subgroup of a hyperbolic group {\em dynamical
malnormality} is equivalent to {\em almost malnormality} which means for all
$g\in G\setminus H$, $H\cap gHg^{-1}$ is finite (see \cite{short}).
However, this is no longer true if we take a relatively quasi convex subgroup of a relatively
hyperbolic group.

\begin{cor}
Suppose $G$ is a relatively hyperbolic group and $H$ is a relatively quasi-convex, dynamically
malnormal subgroup of $G$. Then $H$ has bounded packing in $G$.

In particular, $H$ has finite width in $G$.
\end{cor}

$Proof:$ One just needs to check the conditions of the above theorem.
A relatively hyperbolic group $G$ acts on its Bowditch boundary $\partial G$ which is
a compact space. Now for any relatively quasi-convex subgroup $H$
we can take $A$ to be its limit set. Condition $(1)$ is true for any convergence action
on a compact set. We know that $G$ acts on $\partial G$ as a convergence group and hence
so does $H$. Condition $(2)$ is proved by Dahmani (see Proposition $1.8$ of \cite{dah-conv}).
Condition $(3)$ holds by definition of dynamical malnormality. $\Box$

\begin{cor}
In a relatively hyperbolic group maximal parabolic subgroups have bounded packing.
\end{cor}

This follows easily since we know that the maximal parabolic subgroups are dynamically malnormal
and their limit sets are singleton sets.

\begin{rem}
If we use the fact that quasi-convex subgroup of a hyperbolic group has finite height as proved
in \cite{GMRS} and the Lemma 3.3 of \cite{mahan-CC} then using the proof of the
above theorem we get an alternative proof of the fact that a quasi-convex subgroup of a hyperbolic
group has bounded packing. However, this method fails to deal with the relatively hyperbolic situation.
\end{rem}

\section{Bounded packing of cyclic subgroups generated by rank-1 elements}

This section contains two proofs of the following theorem using completely different sets of ideas.
The first one is an application of Theorem \ref{th1}, whereas
the second proof uses results from Dahmani-Guirardel-Osin(\cite{dah-gui-osin}) and Sisto(\cite{sisto}).

\begin{theorem}({\em Bounded packing of rank-1 cyclic subgroups})\label{th2}
Suppose $G$ is a CAT(0) group and $H$ is a cyclic subgroup generated by a rank-$1$ element. Then $H$ has
bounded packing in $G$.
\end{theorem}

$Proof:$ Let $Y$ be a CAT(0) space on which $G$ acts properly and co-compactly. Let $X=\partial Y$ be the visual boundary of $Y$
with the visual topology. Then $G$ acts on $X$ by homeomorphism. We will check that the pair $(G,H)$ satisfies all
the conditions of Theorem \ref{th1} to finish the proof. For this purpose we shall use the results of Hamenstadt(\cite{hamen-rank1}). Let $A=\{\xi_1,\xi_2\}\subset X$ - the set containing the two fixed points
in  $X$ of the non-trivial isometries in $H$. 

The condition $(1)$ of Theorem \ref{th1} is verified by the fact that rank-$1$ isometries have north-south dynamics on $X$ (Lemma 4.4 of \cite{hamen-rank1}). 

For condition $(2)$ suppose $\{g_nH\}$ is a sequence of distinct cosets of $H$ in $G$ and suppose $\gamma=[\xi_1, \xi_2]$ 
is the geodesic line in $Y$ invariant under $H$. Suppose now that $g_n\xi_i$ converges to $\alpha_i$
for $i=1,2$ and $\alpha_1 \neq \alpha_2$. Fix a point $p\in Y$ and a drop a perpendicular $[p,p_n]$ on to $g_n\gamma$.
Then by Lemma 3.2(2) of \cite{hamen-rank1} the points $p_n$ are uniformly close to the geodesics $[p,g_n\xi_1]$ and 
$[p,g_n\xi_2]$. Since $g_n\xi_i$ converges to $\alpha_i$ and $\alpha_1 \neq \alpha_2$, we have that 
$\angle_p(g_n\xi_1,g_n\xi_2)\geq \epsilon$ for some $\epsilon>0$ for all $n$. Now choose points $q_n\in [p,g_n\xi_1]$,
$r_n\in [p,g_n\xi_2]$
such that $q_n,r_n$ are both uniformly close to $p_n$. Then $d(q_n,r_n)$ is uniformly small but $d(p,q_n),d(p,r_n)$ are
going to infinity. This is so because $\{g_nH\}$ are distinct cosets and hence their distances from $1\in G$ are
going to infinity. This gives a contradiction because suppose $\triangle A_nB_nC_n$ is a comparison triangle in
$\mathbb E^2$ of the geodesic triangle $\triangle pq_nr_n \subset Y$ where $p$ is mapped to $A_n$. Then 
$\angle_{A_n}(B_n, C_n) \geq \angle_{p}(q_n, r_n)\geq \epsilon$ whereas $d(B_n,C_n)$ is uniformly small and $d(A_n,B_n)$
and $d(A_n,C_n)$ are arbitrarily large. Such a sequence of triangles can not exist in $\mathbb E^2$. $\Box$

For condition $(3)$ is a direct consequence of Lemma 4.5 of Caprace(\cite{cap-lec}) and Proposition 4.3 of \cite{hamen-rank1}. $\Box$

\begin{rem}
We note that in general a {\em purely rank-1} subgroup (i.e. one whose non-identity elements are
all rank-1 isometries) of a CAT(0) group may not have
bounded packing. In fact, as in Example 2.22 of \cite{hruska-wise} one can have subgroups of
CAT(-1) groups without bounded packing.
\end{rem}

{\bf An alternative proof of Theorem \ref{th2}:}

Now we sketch an alternative proof of the above theorem using results of Dahmani, Guirardel and
Osin (\cite{dah-gui-osin}) and Sisto (\cite{sisto}). We refer the reader to \cite{dah-gui-osin}
for terminologies. We will need the following theorems as ingredients of the proof.

\begin{theorem}(See Theorem 1.1 and 1.3 of \cite{sisto})\label{th2.1}
Suppose $G$ is a group acting properly and cocompactly on a CAT(0) space $X$ and 
$g\in G$ acts as a rank-1 isometry on $X$. Then there is a virtually cyclic subgroup $E(g)$ of $G$
containing $g$ such that $E(g)$ is hyperbolically embedded in $G$.
\end{theorem}

\begin{theorem}(Corollary 6.36 of \cite{dah-gui-osin})\label{th2.2}
Let $\{H_{\lambda} \}$, $\lambda \in \Lambda$ be a hyperbolically embedded collection of subgroups of
a group $G$. Then for every $\alpha > 0$ there exists finite subsets 
$F_{\lambda} \subset H_{\lambda}\setminus \{1\}$ such that any collection
$\{N_{\lambda} \}$, $\lambda\in \Lambda$ , where $N_\lambda \unlhd H_{\lambda}$ and 
$N_{\lambda}\cap F_{\lambda} = \emptyset$ for every $\lambda \in \Lambda$, is $\alpha$-rotating.
\end{theorem}

\begin{theorem}(Corollary 5.4 of \cite{dah-gui-osin})\label{th2.3}
Let H be a $200$-rotating subgroup of a group $G$. Then the normal subgroup of
$G$ generated by $H$ is a free product of a (usually infinite) family of conjugates of $H$.
\end{theorem}

We first prove the following:

\begin{theorem}\label{th3}
Suppose $H$ is a residually finite hyperbolically embedded subgroup of a group $G$.
Then $H$ has bounded packing in $G$.
\end{theorem}

We start with a lemma:

\begin{lemma}\label{basic}
In an infinite free product of groups $G_1*G_2*\cdots$ each free factor has bounded packing.
\end{lemma}

{\em Proof of the lemma:} Let us show, without loss of generality, that $G_1$ has bounded packing
in $G_1*G_2*\cdots$. Suppose $G_1$ does not have bounded packing. Then there is $D>0$ and an infinite
collection of sets of cosets of $G_1$ which are pairwise at a distance $D$. We know that
$d(xG_1,yG_1)\leq D$ if and only if $d(1, G_1x^{-1}yG_1)\leq D$ whence the double coset has an
element in $B(1,D)$. Since the metric of the group is proper we conclude that there is a finitely
generated subgroup containing $G_1$ in which $G_1$ fails to have bounded packing. Therefore, we 
are reduced to the case of a finite free product. This case follows easily by looking at the
tree of spaces corresponding to the free product decomposition. 
This completes the proof of the lemma. $\Box$

{\em Proof of Theorem \ref{th3}}:
We apply Theorem \ref{th2.2} with $\alpha=200$. Since $H$ is residually finite there is a normal
subgroup $N$ of finite index in $H$ satisfying the condition of Theorem \ref{th2.2}. Now, by Theorem
\ref{th2.3} the normal subgroup in $G$ generated by $N$ is a free product of copies of conjugates
of $N$. It follows that this normal subgroup has bounded packing in $G$ by Lemma \ref{lemma1}(7),
whereas $N$ has bounded packing in the free product by Lemma \ref{basic}. Hence, by Lemma
\ref{lemma1}(3) $N$ has bounded packing in $G$. Again by Lemma \ref{lemma1}(4) we have that
$H$ has bounded packing in $G$. $\Box$

{\em Proof of Theorem \ref{th2}:} 
Let $H=<g>$. By Theorem \ref{th2.1} $E(g)$ is a virtually cyclic subgroup of $G$ which is
hyperbolically embedded in $G$. Hence $E(g)$ has bounded packing in $G$ by Theorem \ref{th3}.
Since $H$ has finite index in $E(g)$ it now follows from Lemma \ref{lemma1}(4) that $H$
has bounded packing in $G$. $\Box$

\section{ Non-existence of bounded packing}

\begin{prop}\label{free}
Let $\mathbb F_2$ denote the free group on two generators.
Then there are finitely generated subgroups of $\mathbb F_2\times \mathbb F_2$ without the
bounded packing property.
\end{prop}

$Proof:$
Let $G$ be any countable group and let $\Delta$ be the diagonal in $G\times G$.
It is clear that any coset of $\Delta$ in $G\times G$ can be written 
uniquely as $(g,1)\Delta$ for some $g \in G$ and any two cosets $(g,1)D$ and $(gh,1)D$ are close
essentially means that there is $x \in G$ such that $xhx^{-1}$ is close to identity in $G$.

We first search for groups $G$ such that the diagonal $\Delta\leq G\times G$ does not have
bounded packing. Towards that goal, it is enough to furnish a group $G$ which has a sequence
of distinct elements $g_n$ such that for all $m,n$, a conjugate of $(g_m)^{-1}.g_n$ is uniformly
close to identity.  Bartholdi, Cornulier and Kochloukova( see \cite{bart-cornu}) 
constructed these types of examples using generalized wreath
products. By their construction we can have examples of the
form $G=\ZZ\wr_X Q$ where the $Q$-action on $X$ is $2$-transitive and $X$ is an infinite set.
Hence, one can take $g_{x}$ to be
the $1$ in the $x$-th copy of $\ZZ$ in $\bigoplus_{x\in X}\ZZ\subset G$ for all $x\in X$.
It follows that for all $x,y\in X$ $g_x^{-1}g_y\in G$ is a conjugate of $g_{x_0}^{-1}g_{y_0}$ for
some fixed
$x_0,y_0\in X$; this means each $g_x^{-1}g_y\in G$ has a conjugate uniformly close to the identity in $G$.
Therefore, for any such group $G$ the diagonal $\Delta\leq G\times G$ does not have bounded packing.
Note that in these examples the groups $G$ are finitely presented too.

Next, choose such a group $G$. Since the group $G$ is finitely generated (presented) there is a surjective
group homomorphism $\phi: \mathbb F_n\rightarrow G$. 
Now we can use Mihailova's construction. That will give us that 
$(\phi\times\phi)^{-1}(\Delta)\leq \mathbb F_n\times \mathbb F_n$ is finitely generated,
where $\phi\times \phi $ is the naturally induced map $\mathbb F_n\times \mathbb F_n\rightarrow G\times G$
and $\Delta\leq G\times G$ is the diagonal.

Since there are natural embeddings of $\mathbb F_n\times \mathbb F_n$ in
$\mathbb F_2\times \mathbb F_2$, we are done. $\Box$

\smallskip
\begin{rem}
Any finitely presented subgroup $G$ of $\mathbb F_2\times \mathbb F_2$ is virtually a product
$H_1\times H_2$ where $H_i\leq \mathbb F_i$ are finitely generated subgroups. We know that
finitely generated  subgroups of free groups are quasi-convex and by the results of \cite{hruska-wise}
that quasi-convex subgroups of hyperbolic groups have bounded packing. 
Hence $H_i\leq \mathbb F_2$ have bounded packing. Therefore, $G$ has
bounded packing in $\mathbb F_2\times \mathbb F_2$. So one can get only finitely generated
counter-examples like above.
\end{rem}

\begin{prop}
There is a CAT(0) group which has a finitely presented subgroup without bounded packing.
\end{prop}

We shall make use of the following theorems for the proof.

{\em {\bf 1-2-3 Theorem}(Baumslag-Bridson-Miller-Short \cite{bbms1})
Suppose that $1\rightarrow N \rightarrow \Gamma \stackrel{p}{\rightarrow} Q\rightarrow 1$
is exact, and consider the fibre product
$$ P := \{(\gamma_1, \gamma_2) | p(\gamma_1) = p(\gamma_2)\} \subset \Gamma\times \Gamma.$$
If $N$ is finitely generated, $\Gamma$ is finitely presented and $Q$ is of type $F$, then $P$
is finitely presented.}

\begin{theorem}(Wise \cite{wise-rips1})
Let $Q$ be a finitely presented group. Then there exists a group $G$ which is the fundamental group
of a compact negatively curved $2$-complex and a finitely generated normal subgroup $N \leq G$
such that $G/N \simeq Q$.
\end{theorem}

{\em Proof of the proposition:}
The proof of this proposition also follows the same line of arguments as that of
Proposition \ref{free} by using the above theorems. First we find a group $Q$ which satisfies
$F_3$ using \cite{bart-cornu} and then apply Wise's construction to get a CAT(-1)
group $G$ with a surjective map $\phi: G\rightarrow Q$ whose kernel is finitely generated.
Finally we use the {\bf 1-2-3} Theorem of Baumslag-Bridson-Miller-Short to conclude
that $H= (\phi\times \phi)^{-1}(\Delta)\leq G\times G$ is finitely presented where
$\Delta\leq Q\times Q$ -the diagonal subgroup, does not have bounded packing. 
Thus $H\leq G\times G$ will not have bounded packing. Since $G$ is a CAT(-1) group $G\times G$ is
a CAT(0) group. $\Box$

\bibliography{BP2.bib}
\bibliographystyle{amsalpha}

\end{document}